\documentstyle[12pt]{article}
\title{\huge\bf On extreme zeros of classical orthogonal polynomials}

\author{   {\bf Ilia Krasikov}\\
	   Brunel University\\
	   Department of Mathematical Sciences\\
	   Uxbridge UB8 3PH United Kingdom\\
	   e-mail: mastiik@brunel.ac.uk
       }

\date{}

\newcommand{\QED}{\hfill$\Box$}

\newcommand{\be}{\begin{equation}}
\newcommand{\ee}{\end{equation}}

\newtheorem{lemma}{Lemma}
\newtheorem{theorem}{Theorem}
\newtheorem{remark}{Remark}

\mathchardef\inn="3232
\renewcommand{\in}{\mbox{$\,\inn\,$}}

\begin{document}

%------------------------- Title Page ----------------------------------

\maketitle
\vspace*{4ex}

\begin{abstract}

Let $x_1$ and $x_k$ be the least and the largest zeros
of the Laguerre or Jacobi polynomial of degree $k.$
We shall establish sharp inequalities 
of the form $x_1 <A, \, \, x_k >B,$ which are uniform in all
the parameters involved. Together with inequalities
in the opposite direction, recently obtained by
the author, this locates the extreme zeros of classical orthogonal
polynomials with the relative precision, roughly speaking, $O(k^{-2/3}).$
\end{abstract}

%\begin{quote}
%\baselineskip3ex
%\end{quote}
%\vspace{1ex}
\noindent
\hspace{5ex}{\bf Keywords:}
orthogonal polynomials, Laguerre polynomials, Jacobi polynomials, zeros
\noindent
%================= Footnotes =================
\footnote{ 2000 \emph{Mathematics Subject Classification} 33C45}
     
%============== Paging control ===============

\thispagestyle{empty}
\addtocounter{page}{-1}
%\newpage

%\newpage
%\baselineskip3.3ex

%------------------------------------------------------------------------%
%                                                                        %
%         1. INTRODUCTION                                                %
%                                                                        %
%------------------------------------------------------------------------%

\renewcommand{\thesection}{\arabic{section}.}
\section{Introduction}
Study of extreme zeros of the Hermite, Laguerre and Jacobi polynomials
has a long history and most of the classical results are
collected in \cite{szego}. But only recently attention has been shifted
to the case when the parameters may vary with the degree $k$ of a polynomial
\cite{assche99, dette, gavro, ismail, klg, leven, msv}.
Most of these results are of the asymptotic nature (with \cite{ismail} and \cite{leven}
being a remarkable exception)
and hold under certain restrictions on the parameters.
Recently the author obtained the following  explicit uniform bounds \cite{kg} (similar inequalities
for the Laguerre case were given earlier in \cite{klg}).
\begin{theorem}
\label{lagr}
Let $ x_1$ and $x_k$ be the least and the largest zero of the Laguerre polynomial
$L_k^{(\alpha)} (x)$ respectively,
$\alpha >-1.$
Then
\begin{equation}
\label{lm}
 x_1 > V^2+3 V^{4/3} (U^2-V^2)^{-1/3}\, ,
\end{equation}
\begin{equation}
\label{lM}
x_k < U^2-3 U^{4/3} (U^2-V^2)^{-1/3}+2 \, ,
\end{equation}
where $V=\sqrt{k+\alpha+1}-\sqrt{k},$ $U=\sqrt{k+\alpha+1}+\sqrt{k} \, .$
\end{theorem}
\begin{theorem}
\label{jac}
Let $x_1$ and $x_k$ be the least and the largest zero of the Jacobi
polynomial $P_k^{( \alpha , \beta )}(x)$ respectively,
$\alpha \ge \beta > -1.$
Then
\begin{equation}
\label{ozjam}
x_1 > \, A+3 (1-A^2 )^{2/3} \, (2R)^{-1/3} \, ,
\end{equation}
\begin{equation}
\label{ozjacM}
x_k <
B-3 (1-B^2 )^{2/3} \, (2R)^{-1/3} +\frac{4q (s+1)}{(r^2+2s+1)^{3/2}}
\, ,
\end{equation}
where
$$s=\alpha +\beta +1 , \; \;
q=\alpha - \beta , \; \; r=2k+\alpha + \beta +1,
\; \; R=\sqrt{(r^2-q^2+2s+1)(r^2-s^2)} \, ,$$
and
$$A=  - \, \frac{R+q(s+1)}{r^2+2s+1} \, , \; \;
B = \frac{R -q(s+1)}{r^2+2s+1} \, .
$$
\end{theorem}

As the zeros of the Hermite polynomials can be easily expressed through the
zeros of the corresponding Laguerre polynomials we will not consider them in this paper.
\\
Previously known results give, roughly speaking,
$V^2 < x_1 <x_k <U^2,$ \cite{ ismail, szego} (see also a survey article \cite{gat})
for Laguerre polynomials,
and $A < x_1 <x_k < B,$ \cite{ismail, leven}
for the Jacobi case. It is also known that these bounds are asymptotically correct
under certain assumptions on the parameters.
On the other hand one can expect that much sharper results similar to these of
Theorems \ref{lagr} and \ref{jac} hold in a more general situation. In particular,
analogous inequality analogous to (\ref{lm}) - (\ref{ozjacM}) are known for the zeros of 
Charlier \cite{kc} and binary Krawtchouk polynomials \cite{k1}.

The aim of this paper is to show that the bounds given by Theorems \ref{lagr} and \ref{jac} 
are essentially sharp, thus locating the extreme zeros of the classical orthogonal polynomials
with a high precision.
Namely we shall establish (in a rather elementary way)
two following theorems giving similar inequalities
in the opposite direction. Our method is based on so-called Bethe ansatz equations,
having some important applications to orthogonal polynomials \cite{ismail1,krasov}.
It is worth also noticing that the above
bounds $V^2 < x_1 <x_k <U^2,$ and $A < x_1 <x_k < B,$ for the Laguerre and Jacobi
polynomials respectively, are an immediate corollary of the Bethe ansatz equation 
we use here (Lemma \ref{thmain} below).
\begin{theorem}
\label{thlag1}
Let $\delta =\frac{1}{k}+\frac{1}{\alpha+1} < \frac{1}{50},$ then
in the notation of Theorem \ref{lagr},
\begin{equation}
\label{eqlagm}
x_1 < V^2+\frac{9V^{4/3}}{ (U^2-V^2)^{1/3} (2-27 \delta ^{2/3})}.
\end{equation}
Let $k \ge 30,$ then
\begin{equation}
\label{eqlagM}
x_k >
U^2-\frac{9U^{4/3} }{2 (U^2-V^2)^{1/3}}
\end{equation}
provided $\alpha \le 2(3+2 \sqrt{3})k-1 ,$
and
\begin{equation}
\label{eqlagM1}
x_k >
U^2-\frac{9U^{4/3} }{(U^2-V^2)^{1/3}(2-3 k^{-2/3})},
\end{equation}
otherwise.
\end{theorem}
\begin{theorem}
\label{jaczer}
Let $ \alpha \ge \beta >-1,$ 
then in the notation of Theorem \ref{jac},
for $ k \ge 5,$
\begin{equation}
\label{jacm}
x_1<
A+9 (1-A^2)^{2/3}(2 R)^{-1/3} \; ,
\end{equation}
and for $k \ge 56,$
\begin{equation}
\label{jacM}
x_k > B- 9 (1-B^2)^{2/3}(2 R)^{-1/3} .
\end{equation}
\end{theorem}

It seems that the bounds in this direction received much less attention.
We will use here some rather weak classical inequalities (\cite{szego}, sec.6.2).

Theorems \ref{lagr} - \ref{jaczer} yield the asymptotics for the extreme zeros
given in the next theorem (in the Jacobi case $x_k$ and $B$ may vanish what leads to
more complicate expressions). The meaning of $O$-terms here is that for sufficiently large $k,$
say $k >100,$ one can replace them by absolute constants.
\begin{theorem}
\label{assym}
(i) In the notation of  Theorem \ref{lagr}, for  sufficiently large $k$
and $\alpha >50,$
the extreme zeros of the Laguerre polynomial 
$L_k^{(\alpha)} (x)$ satisfy
\begin{equation}
\label{aslag1}
\frac{x_1 }{V^2}= 1+O \left((\alpha +1)^{-1/2}(\frac{1}{\alpha +1}+\frac{1}{k})^{1/6} \right),
\end{equation}
\begin{equation}
\label{aslagk}
\frac{x_k }{U^2} = 1-O \left(  k^{-1/6}(k+\alpha )^{-1/2} \right).
\end{equation}
(ii) In the notation of  Theorem \ref{jac}, for  sufficiently large $k$
and $\alpha \ge \beta >-1,$
the extreme zeros of the Jacobi polynomial $P_k^{( \alpha , \beta )}(x)$ satisfy
\begin{equation}
\label{asjac10}
\frac{x_1 }{A}= 1+
O  \left( \left( \frac{(\beta +1)^2}{k(k+ \alpha )(k+ \beta )}\right)^{2/3}\right) ; \; \; r^2 \ge q^2+r^2 ,
\end{equation}
\begin{equation}
\label{asjac11}
\frac{x_1 }{A}=1+
O  \left(\frac{(\beta +1)^{4/3}}{k^{2/3}
(k+ \beta )^{5/6} \sqrt{k+ \alpha}}\right); \; \; r^2 < q^2+s^2 ,
\end{equation}
Let $r^2=q^2+s^2+\gamma (s+1)^{2/3} (r^2-s^2)^{1/3},$ then
\begin{equation}
\label{asjac20}
\frac{x_k}{B} = 1-O \left(\gamma^{-1}+\gamma^{-2/3} k^{-2/9} \right), \, \, 
\, \, \gamma >0;
\end{equation}
\begin{equation}
\label{asjac21}
\frac{x_k}{B} = 1- O \left(( \alpha k )^{-1/3} \right),
\, \, \, \, \gamma < - \; \frac{3(s+1)^{4/3}}{4(r^2-s^2)^{1/3}};
\end{equation}
\begin{equation}
\label{asjac22}
\frac{x_k}{B} = 1- O \left(| \gamma |^{-1}+ | \gamma |^{-1/2} k^{-1/3} \right),
 \;\; \; \; \; - \; \frac{3(s+1)^{4/3}}{4(r^2-s^2)^{1/3}} \le \gamma <0;
\end{equation}
\begin{equation}
\label{asjac23}
|x_k | =O\left( \frac{1}{k^{1/6} \sqrt{k+ \alpha } }\right), \; \; \; \; \; | \gamma | \le 1.
\end{equation}
\end{theorem}
It is worth to compare the obtained inequalities with the classical
results for the fixed values of the parameters.
In particular, in the Laguerre case one has (\cite{szego}, Theorem 6.32, see also \cite{mtn}
for a far-reaching generalization)
$$x_k < \left( \sqrt{4k+2 \alpha +2} -6^{-1/3} (4k+2 \alpha +2)^{-1/6} i_{11} \right)^2,$$
where $6^{-1/3} i_{11}= 1.85575...,$ and $i_{11}$ stands for the least positive zero of
the Airy function.
One can check that for a fixed $\alpha$ 
this differs from (\ref{lM}) only by the better factor $c= 2 \cdot 6^{-1/3} i_{11},$
instead of $3,$ before the second terms of (\ref{lM}).
It is tempting to conjecture that asymptotically
for $k \rightarrow \infty ,$ and uniformly in all the parameters involved,
one should get
the same constant $c$ instead of $3$ before the second terms in all the expressions
(\ref{lm})-(\ref{ozjacM}).

The paper is organized as follows. In the next section we establish rather general
inequalities being our main tool in the sequel. In sections 3 and 4 we will
prove Theorems \ref{thlag1} and \ref{jaczer}, dealing with
Laguerre and Jacobi polynomials respectively. Section 4 also contains a proof of Theorem \ref{assym}.
%%%%%%%%%%%%%%%%%%%%%%%%%%%%%%%%%%%%%%%%%%%%%%%%%%%%%%%%%%%%%%%%%%%%%%%%%%%%%%%%%%%%%%%%%%%%%%%%
\section{Bethe Ansatz Inequalities}
In this section we will consider real polynomials
$f=f(x)$  with only real simple zeros
$x_1 <x_2 <...<x_k ,$ 
satisfying a differential equation 
\begin{equation}
\label{difeq}
f''-2a f'+b f=0
\end{equation}
We suppose here that $a=a(x)$ and $b=b(x)$ are meromorphic functions and none of $x_i$ coincides
with singularities of $a$ or $b.$
For such an $f$ we define the discriminant $\Delta (x)=b(x)-a^2(x),$
and consider the second negative moments of $f$ at its zeros
$$
S(f,x_i) = \sum_{j \ne i} \frac{1}{(x_i-x_j)^2} .
$$
\begin{lemma}
\label{thmain}
\begin{equation}
\label{eq0}
S(f,x_i) = \sum_{j \ne i} \frac{1}{(x_i-x_j)^2}=\frac{ \Delta (x_i )-2 a'(x_i)}{3}.
\end{equation}
\end{lemma}
\begin{proof}
Using the logarithmic derivative and (\ref{difeq})  we get
\begin{equation}
\label{eq5}
\sum \frac{1}{(x-x_j)^2} =- \left( \frac{f'}{f} \right)' =\frac{f'^2-f f''}{f^2}=
\frac{f'^2-2a f' f+b f^2}{f^2}.
\end{equation}
Thus 
$$
S(f,x_i) = \lim_{x \rightarrow x_i} \left(\frac{f'^2-2a f' f+b f^2}{f^2} -\frac{1}{(x-x_i)^2} 
\right).
$$
The result follows on applying four times L'H\^{o}pital's rule
and substituting $f''$ from (\ref{difeq}) at each step.
\QED
\end{proof}
\begin{remark}
Results of this type are called Bethe ansatz equations and are known (or can be routinely
established) in a more general situation and weaker smoothness assumptions. 
We refer to \cite{amm, ismail1, krasov} and the references
therein for a more detailed discussion.
\end{remark}
\begin{lemma}
\label{lem1}
\begin{equation}
\label{ineq}
D (f,x_i,x) = 1+(x-x_i)^2 \left( \frac{ \Delta (x_i )-2 a'(x_i)}{3} - \Delta (x) \right)  >0,
\end{equation}
provided $x \notin [ x_1 ,x_k ].$
In particular, if $a'(x_i) \ge 0,$ then
\begin{equation}
\label{ineq1}
3-2 (x-x_i)^2 \Delta (x_i)+3(x-x_i)^2 \left(\Delta (x_i )-\Delta (x)  \right) >0.
\end{equation}
\end{lemma}
\begin{proof}
From (\ref{eq5}) we have
$$
\frac{1}{(x-x_i)^2} +\sum_{j \ne i} \frac{1}{(x-x_j)^2} = 
(\frac{f'(x)}{f(x)} -a(x))^2 +b(x) -a^2(x) \ge \Delta (x) .$$
Since
$$\sum_{j \ne i} \frac{1}{(x-x_j)^2} < \sum_{j \ne i} \frac{1}{(x_i-x_j)^2} =S(f, x_i ),$$
for $x \notin [ x_1 ,x_k ],$ we obtain
$$\frac{1}{(x-x_i)^2} +S(f, x_i ) > \Delta (x),$$
and (\ref{ineq}),(\ref{ineq1}) follow by Lemma \ref{thmain}.
\QED
\end{proof}
\begin{remark}
Similar arguments can be apply to $x \in [x_1,x_k ],$ say $x_i <x <x_{i+1},$ giving an upper bound
on $x_{i+1}-x_i.$ Indeed,
$$ \Delta (x) \le \frac{1}{(x-x_i)^2}+\frac{1}{(x-x_{i+1})^2} +\sum_{j<i} \frac{1}{(x-x_j)^2}+
\sum_{j>i+1} \frac{1}{(x-x_j)^2} <
$$
$$
\frac{1}{(x-x_i)^2}+\frac{1}{(x-x_{i+1})^2} + \sum_{j < i} \frac{1}{(x_j-x_i)^2}+
\sum_{j > i+1} \frac{1}{(x_j-x_{i+1})^2} <
$$
$$
\frac{1}{(x-x_i)^2}+\frac{1}{(x-x_{i+1})^2} -\frac{2}{(x_{i+1}-x_i)^2}+S_2 (f,x_i )+S_2 (f,x_{i+1} ).
$$
By substituting here $x= \frac{x_i+x_{i+1}}{2},$ one obtains
$$(x_{i+1}-x_i)^2 < \frac{18}{3 \Delta ( \frac{x_i+x_{i+1}}{2})-
\Delta (x_i )-\Delta (x_{i+1})+2 a'(x_i)+2 a'(x_{i+1})},$$
provided the denominator is positive.
\end{remark}
We will solve inequality (\ref{ineq}) for the Laguerre and Jacobi polynomials
in the next section. This will require rather involved calculations
but the following simple heuristic arguments show what type of
bounds may be expected.
\\
Suppose that $\Delta (x)$ has only two real zeros $ y_1 <  y_2 .$
Neglecting the term $2 a'(x),$ we obtain that all the zeros of $f$ are in the interval
$( y_1,  y_2).$ Let $ x_k$ be, say, the largest zero of $f,$ we put
$x_k= y_2 -\epsilon ,$ and choose $x=  y_2 - \frac{5 \epsilon }{9}.$
Now, on omitting higher derivatives of $\Delta ,$ that is putting
$\Delta ( y_2 - \delta )=\Delta ( y_2 )-\delta
\Delta' (  y_2 )= -\delta \Delta' (  y_2 ),$
(\ref{ineq}) can be rewritten as
$$0 < 1+ \frac{16 \epsilon^2}{81} \left( \frac{\Delta (x_k)}{3} - \Delta (x) \right) \approx 
1+\frac{32 \epsilon^3 \Delta' (  y_2 )}{729}.$$
Thus we obtain $x_k >  y_2 +\frac{9}{2} ( 4 \Delta' (  y_2 ) )^{-1/3}.$
Notice that similar heuristic considerations given in \cite{kg} yield in the opposite direction
$x_k <  y_2 +3 ( 4 \Delta' (  y_2 ))^{-1/3},$ (
$\Delta' (  y_2 )$ is negative as $\Delta (  y_2 )=0$).
%%%%%%%%%%%%%%%%%%%%%%%%%%%%%%%%%%%%%%%%%%%%%%%%%%%%%%%%%%%%%%%%%%%%%%%%%%%%%%%%%%%%%%%%%%%%%%%%%%%%%%%%%%

\section{Laguerre Polynomials} 
The Laguerre Polynomials $L_k^{(\alpha)}(x)$ are polynomials orthogonal on $[0, \infty )$ for $\alpha >-1,$
with respect to the weight function $x^{\alpha } e^{-x}.$
The corresponding ODE is
$$
u''-(1-( \alpha+1) x^{-1})u'+k x^{-1} u=0, \; \;\; u=L_k^{(\alpha)}(x).
$$
We also need the explicit representation
\begin{equation}
\label{expll}
L_k^{(\alpha)}(x)= \sum_{i-0}^k {k+ \alpha \choose k-i} \frac{(-x)^i}{i!}.
\end{equation}
Using the notation of Theorem \ref{lagr} we get
$k=\frac{(U-V)^2}{4}, \, \, \alpha =V U-1 ,$ and the condition $\alpha > -1,$
means $ V >0.$

We have 
$a(x)= \frac{x - V U}{2 x} , \, \, a'(x)= \frac{V U}{2 x^2} >0, $
and also
\begin{equation}
\label{eq3}
\Delta (x)= \frac{(U^2-x)(x-V^2)}{4x^2},
\end{equation}

Let $x_1$ and $x_k$ be the least and the largest zeros of $L_k^{(\alpha)}(x)$
respectively. We need the following (rather weak) bounds on $x_1$ and $x_k ,$
(\cite{szego}, sec. 6.2).
\begin{equation}
\label{szlag}
x_1  \le \frac{(\alpha+1)(\alpha +3)}{2k+ \alpha +1}= 
\frac{2 V U (V U+2)}{V^2+U^2}. 
\end{equation}
By (\ref{expll})  we have 
$\sum_{i=0}^k x_i=k(k+ \alpha ),$ implying
$x_1 < k+ \alpha = \frac{(U+V)^2}{4} <x_k.$ 
%Since $S_1(L_k^{(\alpha)},x_1) <0$ it follows $ x_1 <V U.$
Moreover,  as $0 < S(L_k^{(\alpha)},x_i) < \Delta (x_i ) ,$ we get
that all the zeros satisfy $V^2 < x_i <U^2,$  
and thus
\begin{equation}
\label{eqlo}
V^2 <x_1 <  \frac{(\alpha+1)(\alpha +3)}{2k+ \alpha +1} <x_k <U^2 .
\end{equation}
\begin{lemma}
For $V^2 <x < x_1 ,$
\label{lemdvn}
\begin{equation}
\label{eqdn}
\Delta (x_1 )- \Delta (x) < \frac{U^2-V^2}{4 V^4 } (x_1-x).
\end{equation}
For $ x_k <x <U^2 ,$
\begin{equation}
\label{eqdv}
\Delta (x_k )- \Delta (x) < \frac{U^2-V^2}{4 x_k^2} \, (x- x_k).
\end{equation}
\end{lemma}
\begin{proof}
Using that $\frac{(V^2+U^2)x y -V^2 U^2 (x+y)}{x y} ,$ is an
increasing function in $x$ and $y$
we obtain
$$
\frac{\Delta (x_1 )- \Delta (x)}{x_1 -x} = 
\frac{ V^2 U^2 (x+x_1)-(V^2+U^2)x \; x_1}{4 x^2  x_1^2} <
\frac{U^2-V^2}{4 x  x_1} < \frac{U^2-V^2}{4 V^4 };
$$
$$
\frac{\Delta (x_k )- \Delta (x)}{x-x_k} = \frac{(V^2+U^2)x 
x_k-V^2 U^2 (x+x_k)}{4 x^2  x_k^2} 
< \frac{U^2-V^2}{4 x  x_k} <  \frac{U^2-V^2}{4 x_k^2} .
$$
and the result follows.
\QED
\end{proof}

\noindent
$Proof$ $of$ $Theorem$ \ref{thlag1}.
\\
(i) We choose $x=x_1 -\epsilon,$ where $\epsilon = \frac{2 V^{4/3}}{(U^2-V^2)^{1/3}}.$
Then (\ref{ineq1}) and (\ref{lemdvn}) gives
$$0< 3- \frac{\epsilon^2 (U^2-x_1)(x_1-V^2)}{2 x_1^2}+\frac{3 \epsilon^3 (U^2-V^2)}{4 V^4}=
9-\frac{\epsilon^2 (U^2-x_1)(x_1-V^2)}{2 x_1^2} := F(x_1) .$$
We claim that under our assumptions
$F(x)$ has two zeros $y_1 <y_2,$ and 
$x_1 <y_1.$
As $x_1 <x_0=\frac{(\alpha+1)(\alpha +3)}{2k+ \alpha +1},$
it is enough to show that $F(x_0) <0.$ 
Putting $b= \alpha+1,$ we have 
$$F(x_0)= 9+ \frac{2 \epsilon^2}{(b +2)^2}
+\frac{8 \epsilon^2 k (k+b)}{b^2 (b+2)^2}-
\frac{2 \epsilon^2 k (k+b )}{b^2}.
$$
Here
$$
\frac{2 \epsilon^2}{(b +2)^2}+\frac{8 \epsilon^2 k (k+b)}{b^2 (b+2)^2} <
\frac{8 \epsilon^2 (k+b)^2}{b^4} < \frac{16 }{b} \delta^{1/3} <
16 \delta^{4/3},  $$
and
$$\frac{2 \epsilon^2 k (k+b )}{b^2} = \left( \frac{16 b k (b+k)}{U^4} \right)^{2/3} >
\left( \frac{b k }{k+b} \right)^{2/3} = \delta^{-2/3}.  $$
Now it is left to check that  
$9+ 16 \delta^{4/3} - \delta^{-2/3} <0,$
for $ \delta < \frac{1}{50},$
proving the claim.
For $y_1$ we get
$$y_1=V^2+\frac{9V^2}{h(1+\sqrt{1-18 V^2 U^2 h^{-4}} \;)-9},$$ 
where $h=V^{2/3}(U^2-V^2)^{1/3},$
and
$$
h(1+\sqrt{1-18 V^2 U^2 h^{-4}} \; )-9 >2h-\frac{18V^2 U^2}{h^3}-9 < 2h-\frac{27 U^2}{U^2-V^2}=
$$
$$
2h \left(1-\frac{27 U^{8/3}}{2 b^{2/3} (U^2-V^2)^{4/3}} \right) <h (2-27 \delta^{2/3}  ). 
$$
As 
$2-27 \cdot 50^{-2/3} >0,$
the result follows.
\\
(ii)
We choose $x=x_k -\epsilon,$ where $\epsilon = \frac{2 U^{4/3}}{(U^2-V^2)^{1/3}}.$
By (\ref{ineq1}) and (\ref{lemdvn}) we have
$$0 <3- \frac{\epsilon^2 (U^2-x_k )(x_k-V^2)}{2 x_k^2}+\frac{3 \epsilon^3 (U^2-V^2)}{4 x_k^4}=
3+\frac{6 U^4}{x_k^2}-\frac{\epsilon^2 (U^2-x_k)(x_k-V^2)}{2 x_k^2}<$$
$$\frac{9 U^4}{x_k^2}-\frac{\epsilon^2 (U^2-x_k)(x_k-V^2)}{2 x_k^2}.
$$
Thus
$$F(x_k):=18x_k^2-\epsilon^2 (U^2-x_k)(x_k-V^2) >0.$$
The equation $F(x)=0,$ has two zeros,  $y_1 <y_2,$ and $x_k >y_2.$
Indeed, as $x_k >x_0=\frac{(V+U)^2}{4},$ it is enough to check
$F(x_0) <0.$ We have
$$4 F(x_0) =72U^4-\frac{\epsilon^2 (U-V)^2 (3 U^2+10 V U+3 V^2)}{4} \le
72U^4-3 U^{8/3}(U^2-V^2 )^{4/3} \le $$
$$
3 U^{8/3}(U+V)^{4/3}( 24-(U-V)^{4/3}) =3 U^{8/3}(U+V)^{4/3}( 24-(4k)^{2/3}) <0,
$$
for $k \ge 30.$
Thus, 
$$x_k >y_2=U^2-\frac{9 U^2}{9+U^{2/3}(U^2-V^2)^{1/3} \left(
1+\sqrt{1-18V^2 U^{-2/3}(V^2+U^2)^{-4/3}}\right)} >$$
$$U^2-\frac{9U^2 }{2 U^{2/3}(U^2-V^2)^{1/3}+9- \frac{18V^2}{U^2-V^2}}.$$
Finally, $9- \frac{18V^2}{U^2-V^2} \ge 0,$ if  $\alpha \le 2(3+2 \sqrt{3})k-1 ,$
proving (\ref{eqlagM}). Otherwise,
$$
2 U^{2/3}(U^2-V^2)^{1/3}- \frac{18V^2}{U^2-V^2} =
U^{2/3}(U^2-V^2)^{1/3} \left( 2-\frac{18V^2}{U^{2/3}(U^2-V^2)^{4/3}} \right) > $$
$$U^{2/3}(U^2-V^2)^{1/3}
\left( 2-\frac{9b^2}{2^{5/3}k^{2/3}(k+b)^2 } \right) > 
U^{2/3}(U^2-V^2)^{1/3}(2-3 k^{-2/3}), $$
and (\ref{eqlagM1}) follows.
\QED
%%%%%%%%%%%%%%%%%%%%%%%%%%%%%%%%%%%%%%%%%%%%%%%%%%%%%%%%%%%%%%%%%%%%%%%%%%%%%%%%%%%%%%%%%%%%%%%%%%%%%%%%%
\section{Jacobi Polynomials} 
The Jacobi Polynomials $P_k^{( \alpha , \beta )}(x)$ are polynomials orthogonal on $[-1,1]$
for $\alpha , \beta >-1,$
with respect to the weight function $(1-x)^{\alpha}(1+x)^{\beta}.$
The corresponding ODE is
$$
u''- \frac{( \alpha +\beta +2)x+\alpha-\beta}{1-x^2} u'+\frac{k(k+\alpha + \beta +1)}{1-x^2} u=0, \; \; \; u=P_
k^{( \alpha , \beta )}(x).  $$
We will use the notation of Theorem \ref{jac} and put $p=r^2+2s+1$ throughout this section.
We have
\begin{equation}
\label{eq4}
\Delta (x)= - \, \frac{p x^2+2q(s+1)x+s^2+q^2-r^2}{4(1-x^2)^2}=
\frac{p(x-A)(B-x)}{4(1-x^2)^2}.
\end{equation}
As
$$a'(x)= \frac{((\alpha+ \beta +2)x+\alpha - \beta )^2+
4(\alpha+1)(\beta+1)}{2(\alpha +\beta +2)(1-x^2)^2} >0,$$
we can use (\ref{ineq1}) and moreover,
as $\Delta ( x_i) >0,$ we obtain
\begin{equation}
\label{eqhja}
A    < x_i < B 
\end{equation}
In the opposite direction it is known (\cite{szego}, sec. 6.2)
\begin{equation}
\label{szjac}
x_1 <- \; \frac{2k + \alpha - \beta -2}{2k+ \alpha + \beta}  < 
\frac{2k +\beta - \alpha -2}{2k+ \alpha+\beta} < x_k .
\end{equation}
It is easy to show (see e.g. \cite{kg}) that $ x_1 <0 , $
for $\alpha \ge \beta .$
\begin{lemma}
\label{jacraz}
For $A<x <x_1,$
$$\Delta (x_1)- \Delta (x)  < 
\frac{R (x_1 -x)}{2(1- x_1^2 )^2}. $$
For $x_k <x <B,$
$$\Delta (x_k)- \Delta (x)<
\frac{R (x -x_k)}{2(1- x_k^2 )^2}.$$
\end{lemma}
\begin{proof}
We have
$$\Delta (x_1)- \Delta (x)  < \frac{p}{4(1- x_1^2 )^2} \; \left( 
(x_1-A)(B-x_1) -(x-A)(B-x)\right) <$$ 
$$ \frac{p}{4(1- x_1^2 )^2} \; (B-A)(x_1 -x)=
\frac{R (x_1 -x)}{2(1- x_1^2 )^2}.  $$
$$\Delta (x_k)- \Delta (x) 
< \frac{p}{4(1- x_k^2 )^2} \; \left( 
(x_k-A)(B-x_k) -(x-A)(B-x)\right) <$$
$$
\frac{p}{4(1- x_k^2 )^2} \; (B-A)(x -x_k)=
\frac{R (x -x_k)}{2(1- x_k^2 )^2}.
$$
\QED
\end{proof}

\noindent
$Proof$ $of$ $Theorem$ \ref{jaczer}.
\\
(i)  Choose $\epsilon= \frac{(2-2 A^2 )^{2/3}}{R^{1/3}},$
and put $x= x_1- \epsilon .$ Then $x> A,$ 
otherwise there is nothing to prove.
Using the previous lemma and (\ref{ineq1}) we obtain
$$
0 < 3- \frac{ \epsilon^2 p(x_1- A)(
B-x_1)}{2(1-x_1^2 )^2}+ \frac{3 \epsilon^3 R}{2(1- x_1^2 )^2} <
$$
$$
\frac{3 (1-A^2 )^2}{(1-x_1^2 )^2}- \frac{ \epsilon^2 p(x_1- A)(
B-x_1)}{2(1-x_1^2 )^2}+ \frac{3 \epsilon^3 R}{2(1- x_1^2 )^2} \; .
$$
Thus, we get
\begin{equation}
\label{uravj}
18(1-A^2 )^2-\epsilon^2 p (B -x_1)(x_1- A) := F(x_1) >0.
\end{equation}
We shall show that this quadratic has two real zeros $z_1<z_2 , $ and $x_1<z_1.$
For, it is enough to prove
$ F(\frac{A+B}{2})<0,$
and $\frac{A+B}{2}.$ The last claim follows from (\ref{szjac}), as 
$$
x_1 < - \; \frac{2k+ \alpha- \beta-2}{2k+ \alpha+ \beta} <\frac{A+B}{2} .$$
Indeed, $\alpha , \beta >-1,$ and we obtain
$$
\frac{A+B}{2} -x_1 >
\frac{A+B}{2}+ \frac{2k+ \alpha- \beta-2}{2k+ \alpha+ \beta}=
$$
$$
\frac{4(2k^3+(3 \alpha+ \beta +4)k^2+(\alpha^2+ \alpha \beta+
4 \alpha+4 \beta+4)k+(\alpha+1)(\alpha+\beta+2))}{(r-1)p}>0.
$$
Now we have
$$
F(\frac{A+B}{2})=72 (1-A^2 )^2-\epsilon^2 p (B -A)^2,
$$
and it is negative whenever
\begin{equation}
\label{prom1}
2 R^4 > 729 p^3 (1-A^2 )^2 .
\end{equation}
As 
$$
\frac{d}{d q} \left( \frac{1-A^2 }{R^2} \right) =
- \; \frac{2 ((s+1)^2-q^2)}{(p-q^2)R(q R+(s+1)(p-q^2))} <0,
$$
and for $q=0,$
$$
\frac{1-A^2 }{R^2}= \frac{(s+1)^2}{p^2 (r^2-s^2)}.
$$
We have
$$
\frac{p^3 (1-A^2 )^2}{ R^{4}} <\frac{(s+1)^4}{p (r^2-s^2)^2} < 
\frac{(s+1)^4}{16 k^2(k+s)^2(2k+s)^2}<
\frac{1}{16k^2} <\frac{2}{729},
$$
provided $k \ge 5.$
This proves (\ref{prom1}) and, thus, $x_1 <z_1.$
\\
Finally, solving $F(x)=0,$ we obtain
$$
x_1 < A+ \frac{18 (1-A^2 )^2}{\epsilon^2 R \left(
1+\sqrt{1-\frac{18p (1-A^2 )^2}{ \epsilon^2 R^2 }} \; \right)} <A+ 
\frac{18 (1-A^2 )^2}{\epsilon^2 R } =
A+\frac{9 (1-A^2)^{2/3} }{(2 R)^{1/3}}.
$$
\\
(ii)
Choose $\epsilon= \frac{(2-2 B^2 )^{2/3}}{R^{1/3}},$
and put $x= x_k+ \epsilon .$
Similarly to the previous case we get
\begin{equation}
\label{ineqM2}
18(1-B^2 )^2-\epsilon^2 p (B-x_k)(x_k- A) :=F(x_k)>0.
\end{equation}
We shall show that $x_k$ is greater than the largest zero of $F(x)=0.$
To prove this we establish $F(x_0)<0,$
where
$x_0 = \frac{2k+\beta - \alpha -2}{2k+\alpha+\beta } <x_k,$ 
by Lemma \ref{szjac}.
For 
it is enough to show
$$
G= \left( \frac{18 (1-B^2)^2}{ \epsilon^2 p (B-x_k)(x_k- A) } \right)^3 =
\frac{729 R^2 (1-B^2)^2}{2 (p (B-x_k)(x_k- A) )^3} <1.
$$
We have
$$
\frac{d}{d q} \left( \frac{1-B^2 }{R^2} \right) =
\frac{2 (r^2-s^2)(q(r^2-s^2)+(s+1)R)}{p R^4} > 0.
$$
As $q= \alpha -\beta < \alpha +\beta +2=s+1,$
we obtain
$$
\frac{1-B^2}{R^2} < \frac{4(s+1)^2}{p^2 (r^2-s^2)},
$$
that is 
$$1-B^2 < \frac{4(s+1)^2 R^2}{p^2 (r^2-s^2)}.$$
We also have
$$
p(B-x_0)(x_0 -A)= \frac{16( \alpha +1)((k-1)(\alpha +1)+
k(k+ \beta )(2k+\alpha+\beta))}{(2k+\alpha+\beta)^2} \ge
$$
$$
\frac{16k ( \alpha +1)(k+ \beta )}{2k+\alpha+\beta+1} > \frac{2(r-q-1)(s+q+1)(r-s)}{r} .
$$
Therefore we obtain
$$
G< \frac{729 r^3 (s+1)^4 (r+s)(p-q^2)^3}{p^4(r-q-1)^3 (r-s)^2 (s+q+1)^3} <
\frac{729 (s+1) (r+s)(p-q^2)^3}{r^5 (r-q-1)^3 (r-s)^2}.
$$
The last expression is an increasing function in $q,$ and substituting $q=s+1,$
we get
$$
G< \frac{729 (s+1)(r-s)(r+s)^4}{r^5 (r-s-2)^3} =\frac{2916k(s+1)(k+s)^4}{(k-1)^3(2k+s)^5} 
< \frac{2916k}{(k-1)^3} <1,
$$
for $k \ge 56.$
Finally, solving $F(x_0)=0,$ we obtain
$$x_k > B- \frac{18 (1-B^2 )^2}{\epsilon^2 R \left(
1+\sqrt{1-\frac{18p (1-B^2 )^2}{ \epsilon^2 R^2 }} \; \right)} >B-
\frac{18 (1-B^2 )^2}{\epsilon^2 R } =
B-\frac{9 (1-B^2)^{2/3} }{(2 R)^{1/3}}.
$$
\QED
\begin{remark}
More accurate calculations show that in fact (\ref{jacM}) holds for $k \ge 20,$ instead of $56.$
It is also easy to improve the constant $9$ in (\ref{jacm}), (\ref{jacM}) to $\frac{9}{2-o(1)},$
similarly to the Laguerre case.
\end{remark}
\noindent
$Proof$ $of$ $Theorem$ \ref{assym}.
The asymptotics for the Laguerre case is an easy exercise,
here we will establish (\ref{asjac10})-(\ref{asjac23}). 
\\
Notice that the inequality $r^2  \ge s^2+q^2$ is equivalent to $R \ge q(s+1).$
We also observe that the last term in (\ref{ozjacM}) may be ignored.
Indeed,
$$
1-B^2= \frac{(q+s+1)^2 (R+p-q(s+1))}{p(R+p+q(s+1))}, $$
and this is an increasing function on $R.$
As $q <s+1,$ we get $R >r^2-s^2,$
what implies
$$1-B^2 >\frac{(q+s+1)^2(2r^2-s^2-q (s+1))}{p(2r^2-s^2+q (s+1))}
>\frac{2(\alpha +1 )^2(2r^2-s^2-q(s+1))}{p^2} >
$$
$$
\frac{(\alpha +1 )^2(k +\alpha) (k+ \beta )}{p^2}.
$$
Now calculations yield 
$$\left( \frac{q(s+1)R^{1/3}}{p^{3/2}(1-B^2)^{2/3}}\right)^6 < c \; 
\frac{  q^6 (s+1)^6 k}{(\alpha +1 )^8 (k+ \alpha )^4 (k+\beta )^3},$$
for some positive constant $c$.
This expression is a decreasing function in $\beta$ and for $\beta =-1,$ is 
$O \left(\frac{( \alpha +1 )^4}{k^2 (k+ \alpha )^4} \right) .$
Thus the last term in (\ref{ozjacM}) is negligible whenever $k \rightarrow \infty .$
\\
$Proof$ $of$ (\ref{asjac10}). As $R \ge q(s+1),$
we have
$| A| >  \frac{R}{2 r^2} ,$
and
$$1-A^2 <2( 1+A) = \frac{2 (s+1-q)^2}{R+p-q(s+1)} < \frac{8 (\beta +1)^2}{p} <
\frac{16(\beta +1)^2}{r^2}.$$
Therefore,
$$
\left( \frac{(1-A^2)^{2/3}}{|A| R^{1/3}} \right)^{3/2} < \frac{128 (\beta +1)^2 r} {R^2 }
<\frac{256 (\beta +1)^2 r} {(r^2-q^2)(r^2-s^2)} <\frac{32(\beta +1)^2}{k( k+\alpha )(k+\beta)},
$$
and (\ref{asjac10}) follows.
\\
$Proof$ $of$ (\ref{asjac11}). As $R < q(s+1),$
we get  $q^2 >r^2-s^2 .$
This yields $-1< \beta< 2k+\alpha -2 \sqrt{k(2k+2 \alpha+1) }, $
$\alpha >2k-1+2\sqrt{k(2k-1)},$ and $k <\alpha/2 .$
Thus, $s$ is a large positive number and 
$| A| >  \frac{q s}{ r^2} .$
Now, using $R >r^2-s^2,$ we obtain
$$1-A^2 < \frac{2 (s+1-q)^2}{R+p-q(s+1)}  < \frac{8(\beta +1)^2}{2r^2-s^2-q(s+1)}<
 \frac{4(\beta +1)^2}{\alpha (k+ \beta )}.
$$
This yields
$$
\left( \frac{(1-A^2)^{2/3}}{|A| R^{1/3}} \right)^6 < 
\frac{256 r^{12} (\beta +1)^8}{ \alpha^4 (k+ \beta)^4 q^6 s^6 (r^2-q^2)(r^2-s^2)} <
\frac{10^8 (\beta +1)^8}{k^4 (k+\alpha)^3 (k+ \beta)^5 }, 
$$
and the result follows.
\\
$Proof$ $of$ (\ref{asjac20}).
The condition $q^2=r^2-s^2 - \gamma (s+1)^{2/3}(r^2-s^2)^{1/3}, \, \, \gamma >0,$
implies that $R > q(s+1),$ and $B >0.$
Rewriting $B$ as $\frac{r^2-s^2-q^2}{R+q(s+1)} $ we obtain
$B >\frac{r^2-s^2-q^2}{2R}.$
We also have
$$
1-B^2 <2(1-B)= \frac{8(\alpha +1 )^2}{R+p+q(s+1)} <  \frac{8(\alpha +1)^2}{r^2}.
$$
Hence
$$\left( \frac{(1-B^2)^{2/3}}{B R^{1/3}} \right)^3
=\frac{512 (\alpha +1 )^4 }{r^4} \left( \gamma^{-3}+ \frac{(r^2-s^2)^{1/3}}{\gamma^2 (s+1)^{4/3}}\right) <
$$
$$
512 \gamma^{-3}+ \frac{900 k^{1/3} (\alpha +1 )^4 (k+ \alpha + \beta +1)^{1/3} }{ \gamma^2 (\alpha + \beta +2)^{4/3}
(2k+ \alpha + \beta +1)^4} .
$$
The second term here is a decreasing function in $\beta >-1 ,$
and does not exceed
$$
\frac{900 k^{1/3} (\alpha +1)^{8/3} (k+ \alpha )^{1/3}}{\gamma^2 (2k+ \alpha )^4}< 900 k^{-2/3} \gamma^{-2},
$$
and the result follows.
\\
$Proof$ $of$ (\ref{asjac21}), (\ref{asjac22}).
In those case  $k < \frac{\sqrt{2 \alpha^2+2 \alpha+1}-\alpha}{2},$
$\alpha >2k-1+ 2\sqrt{k(2k-1)},$ and so $\alpha$ is large. Therefore,
$$0 >r^2-q^2-s^2 >r^2-2(s+1)^2 >r^2-4 s^2 ,$$
and hence $ s<r <2s.$
By $q <s+1,$ it follows 
$$- \gamma < \frac{2s^2+2s+1-r^2}{(s+1)^{2/3}(r^2-s^2)^{1/3}} <\frac{(s+1)^{4/3}}{(r^2-s^2)^{1/3}}.$$
Rewriting $B$ as $\frac{r^2-q^2-s^2}{R+q(s+1)} ,$
and using $R <q(s+1), $ we have
$B >\frac{r^2-q^2-s^2}{2q(s+1)} .$
Now,
$$
\left( \frac{(1-B^2)^{2/3}}{B R^{1/3}} \right)^6 \le B^{-6} R^{-2} < 
\frac{64 q^2 (s+1)^6}{R^2 (r^2-q^2-s^2)^6} =
$$
$$
\frac{64 (s+1)^{4/3} \left((r^2-s^2)^{2/3}-\gamma (s+1)^{2/3} \right)^3}{\gamma^6 
(r^2-s^2)^2  \left((s+1)^{4/3}+\gamma  (r^2-s^2)^{1/3}  \right)} <$$
$$
\frac{- 256 (s+1)^{10/3} }{\gamma^3 (r^2-s^2)^2  
\left((s+1)^{4/3}+\gamma  (r^2-s^2)^{1/3}  \right)} +
\frac{256 (s+1)^{4/3} }{\gamma^6 \left((s+1)^{4/3}+\gamma  (r^2-s^2)^{1/3}  \right)} =I_1+I_2.
$$
Now we shall consider two cases corresponding to the restrictions in (\ref{asjac22}) and (\ref{asjac21}).
If $-\gamma \le \frac{3(s+1)^{4/3}}{4(r^2-s^2)^{1/3}},$ that is $q^2 < r^2 -\frac{s^2-6s-3}{4},$
then 
$$I_1 <\frac{- 1024 (s+1)^2}{\gamma^3 (r^2-s^2)^2}
<\frac{- 64}{ \gamma^3 k^2}.$$
Otherwise, using $-\gamma <\frac{2s^2+2s+1-r^2}{(s+1)^{2/3}(r^2-s^2)^{1/3}},$  and
$k <\frac{\alpha}{4},$ for large $\alpha,$ we get
$$I_1 <\frac{256(s+1)^6}{ (r^2-s^2)^2 (2s^2+2s+1-r^2)}<\frac{128(s+1)^6}{ \alpha^6 k^2 (k+s)^2}=
O\left( \alpha^{-2} k^{-2} \right).$$
Similarly,
$ I_2 = O\left( \gamma^{-6} \right), $ if $- \gamma \le \frac{3(s+1)^{4/3}}{4(r^2-s^2)^{1/3}},$
$$ I_2 = O\left(\frac{(r^2-s^2)^2}{(s+1)^8} \right)= O\left(k^2 a^{-6}\right),
$$
if $\frac{3}{4}<  \frac{- \gamma (r^2-s^2)^{1/3}}
{(s+1)^{4/3}}\le \frac{6}{7},$
and
$I_2=\left(k \alpha^{-5} \right),$ otherwise.
These readily yield (\ref{asjac21}), (\ref{asjac22}).
\\
$Proof$ $of$ (\ref{asjac23})
In this case
$$|B|=|\; \frac {r^2-q^2-s^2}{R+q(s+1)} \; | < 
\frac{| \gamma | (s+1)^{1/3} }{(r^2-s^2)^{1/6} \sqrt{(s+1)^{4/3}+\gamma (r^2-s^2)^{1/3}}} =$$
$$
O\left(\frac{|\gamma |}{(s+1)^{1/3} (r^2-s^2)^{1/6} }\right)= 
O\left( \frac{1}{k^{1/6} \sqrt{k+ \alpha } }\right),
$$
and $R= s \sqrt{r^2-s^2} \left( 1+o(1)\right).$
Thus,
$$
(1-B^2)^{2/3} R^{-1/3}= \frac{\left( 1-o(1)\right)}{s^{1/3} 
(r^2-s^2)^{1/6}}=O\left( \frac{1}{k^{1/6} \sqrt{k+ \alpha } }\right),
$$
and (\ref{asjac23}) follows.
\QED
%%%%%%%%%%%%%%%%%%%%%%%%%%%%%%%%%%%%%%%%%%%%%%%%%%%%%%%%%%%%%%%%%%%%%%%%%%%%%%%%%%%%%%%%%%%%%%%%%%%%%%%%%

\end{document}